%% file: root.tex
\documentclass{ifacconf}

\input{Set/packages}

\input{Set/commands}
\usepackage{natbib}        
\usepackage{url}
\begin{document}
\begin{frontmatter}

\title{Koopman-based Control for Stochastic Systems: Application to Enhanced Sampling}

\thanks[footnoteinfo]{Lei Guo is supported by Deutsche Forschungsgemeinschaft through Research Training Group 2297 ``Mathematical Complexity Reduction'' (MathCoRe).}

\author[First,Second]{Lei Guo} 
\author[Third]{Jan Heiland} 
\author[First]{Feliks Nüske}

\address[First]{Max-Planck-Institute for Dynamics of Complex Technical Systems, Magdeburg, Germany}
\address[Second]{Otto von Guericke Universität Magdeburg, Germany}
\address[Third]{Technische Universität Ilmenau, Germany}

\begin{abstract}   
We present a data-driven approach to use the Koopman generator for prediction and optimal control of control-affine stochastic systems. We provide a novel conceptual approach and a proof-of-principle for the determination of optimal control policies which accelerate the simulation of rare events in metastable stochastic systems.
\end{abstract}


\end{frontmatter}

\section{Introduction}
\input{Contents/0_Intro}

\section{Koopman Theory for Control-Affine Systems}

\subsection{Stochastic Systems and Koopman Operators}
\input{Contents/1_Stochastic_Systems_and_Koopman_Operators}

\subsection{Learning Methods} 
\input{Contents/2_Learning_methods}

\subsection{Solution of the Control Problem}
\input{Contents/3_Solution_of_Control_Problem}

\section{Results}
\input{Contents/4_Results}

\subsection{Prediction}\label{sec prediction}
\input{Contents/5_Prediction}

\subsection{Tracking}
\input{Contents/6_Tracking}

\subsection{Loss Functions for Enhanced Sampling}
\input{Contents/7_Loss_Functions}

\subsection{Comparison}
\input{Contents/8_Comparison}

\section{Conclusion}
\input{Contents/9_Conclusion}

\section*{Data Availability Statement}
Codes and data to reproduce the results and figures shown in this manuscript are available from the following public repository:\\
\url{https://doi.org/10.5281/zenodo.13838629}

\bibliography{ifacconf}             
\end{document}

%% file: Set/packages.tex
\usepackage[utf8]{inputenc}
\usepackage{graphicx}
\usepackage{subfigure}
\usepackage[export]{adjustbox}
\usepackage[font=footnotesize,labelfont=bf]{caption}
\usepackage{amsmath}
\usepackage{amssymb}
\usepackage{color}
\usepackage{enumerate}

\usepackage[dvipsnames]{xcolor}
\usepackage{tikz}
\usetikzlibrary{positioning}
\usepackage{tkz-euclide}

\usepackage{tabularx}

\usepackage{algorithm}
\usepackage{algpseudocode}

\usepackage{todonotes}

%% file: Set/commands.tex
\newcommand{\htop}{\mathrm{H}}                                   

\newcommand{\cL}{\mathcal{L}}



%% file: Contents/0_Intro.tex
Stochastic dynamics with metastability are a recurring theme in many scientific disciplines, for instance, in simulations of macro-molecules, in climate systems, and in applications of  uncertainty quantification. Metastability describes the existence of long-lived macro-states in a dynamical system's state space, such that transitions between these macro-states are rare events.

As is well known, metastability can be quantified using the dominant spectrum of the associated Koopman operator~[\cite{davies_metastable_1982,dellnitz_approximation_1999,schutte_direct_1999}]. Numerical methods for spectral approximation of the Koopman operator based on variational principles~[\cite{noe_variational_2013}] and the extended dynamic mode decomposition (EDMD) algorithm~[\cite{williams_data-driven_2015}] have been developed in recent years~[\cite{klus_data-driven_2018}]. On the other hand, metastability is also closely related to control systems. There is a wide range of biased sampling algorithms, which seek to overcome the rare event nature of the dynamics using a time-dependent input, see e.g.~[\cite{mehdi_ml_sampling_2024}] for a recent review. In fact, the exit time from a metastable state can be written as the value function of a specific stochastic optimal control problem~[\cite{hartmann_variational_2017}].

In this study, we join the ideas of Koopman-based modeling and biased sampling. The key ingredient is the generator extended dynamic mode decomposition algorithm (gEDMD)~[\cite{klus_data-driven_2020}], a variant of EDMD to approximate the Koopman generator. It was shown in~[\cite{peitz_data-driven_2020}] that for control-affine stochastic differential equations (SDEs), application of gEDMD reduces the Kolmogorov backward equation into an ODE that is bi-linear in expectation and input. This simplified structure can be utilized for designing controllers which are geared towards accelerated sampling of rare events.

This paper is a first proof-of-principle that the gEDMD method for control-affine SDEs can be used to a) accurately predict the expectation of observable functions of interest for fixed control input; b) solve optimal control problems (OCPs) with integrated running cost and terminal cost; c) design OCPs which enforce accelerated transitions between metastable states. The paper is organized as follows: Section 2 provides a brief introduction to Koopman operators for (controlled) stochastic systems, the gEDMD method, and its application to optimal control problems. Section 3 presents the numerical results for various cost functions in OCPs. The conclusion is located in Section 4.

%% file: Contents/1_Stochastic_Systems_and_Koopman_Operators.tex
\subsubsection{Stochastic Differential Equations (SDEs)}
We start by considering a time-homogeneous stochastic differential equation
\begin{equation}\label{sde}
dX_{t}=b(X_{t})dt+\sigma(X_{t})dW_{t},\,\,\,X_{0}=:x
\end{equation}
with the state $X_{t}\in\mathbb{R}^{n}$,  the drift vector field $b\colon\mathbb{R}^{n}\to\mathbb{R}^{n}$, the diffusion matrix field $\sigma\colon\mathbb{R}^{n}\to\mathbb{R}^{n\times s}$ and the Brownian motion $(W_{t})_{t\ge0}$ in $\mathbb{R}^{s}$. Moreover, to guarantee the global existence of solutions to (\ref{sde}), we assume that $b,\sigma$ satisfy standard Lipschitz properties [Thm.5.2.1 in~\cite{oksendal_stochastic_2003}].

\subsubsection{Koopman Operators}
The stochastic Koopman operator $\mathcal{K}^{t}\colon\mathcal{L}^{\infty}(\mathbb{R}^{n})\to\mathcal{L}^{\infty}(\mathbb{R}^{n})$, associated with an observable $g\colon\mathbb{R}^{n}\to\mathbb{C}$, is defined in, e.g. [\cite{mezic_spectral_2005}] as follows
\begin{equation}\label{operator}
\mathcal{K}^{t}g(x)=\mathbb{E}^{x}(g(X_{t})),
\end{equation}
where $\mathbb{E}^{x}(\cdot)$ denotes the conditional expectation given an initial value $x$. The generator $\mathcal{L}$ with its domain $\mathcal{D}(\mathcal{L})$ of the stochastic Koopman operator $\mathcal{K}^{t}$ is defined by
$$\mathcal{L}\phi:=\lim_{t\to0}\frac{\mathcal{K}^{t}\phi-\phi}{t},\,\,\forall\phi\in\mathcal{D}(\mathcal{L}).
$$
It can be shown by using Ito's lemma (see Section 4.4.2 of ~[\cite{evans_introduction_2012}]) that $\mathcal{L}$ acts as a second-order differential operator, i.e.
\begin{equation}\label{generator}
    \mathcal{L}\phi=b\nabla_{x}\phi+\frac{1}{2}\Sigma\colon\nabla_{x}^{2}\phi,
\end{equation}
where the diffusion matrix is given by $\Sigma:=\sigma\sigma^{T}\in\mathbb{R}^{n\times n}$ and $\colon$ is the Frobenius inner product. Semigroup theory implies that the conditional expectation is the solution of so-called backward Kolmogorov equation, i.e.
\begin{equation}\label{generator pde}
\frac{\partial }{\partial t}\mathbb{E}^{x}(\phi(X_{t}))=\mathcal{L}\mathbb{E}^{x}(\phi(X_{t}))
\end{equation}

with the initial value $\mathbb{E}^{x}(\phi(X_{0}))=\phi(x)$.

\subsubsection{Koopman Generator for Controlled SDEs}
Now, consider a controlled stochastic differential equation with control-affine drift defined by
\begin{equation}\label{control stocha system}
dX_{t}=(b(X_{t})+\sum_{i=1}^{p}G_{i}(X_{t})u_{i}(t))dt+\sigma(X_{t})dW_{t},
\end{equation}
where $u\colon\mathbb{R}_{\ge0}\to\mathbb{R}^{p}$ is the input signal and $
G_{i}\colon\mathbb{R}^{n}\to\mathbb{R}^{n}$, $i=1,\cdots,p$, are Lipschitz continuous vector fields. The Koopman generator for such systems, derived from (\ref{generator}), is given as follows:
\begin{equation}\label{controlled generator}
    \mathcal{L}_{u}\phi=(b+\sum_{i=1}^{p}G_{i}u_{i}(t))\nabla_{x}\phi+\frac{1}{2}\Sigma\colon\nabla_{x}^{2}\phi.
\end{equation}
As in~\eqref{generator pde}, the operator $\cL_u$ evolves the following partial differential equation (see Section III of~\cite{fleming_controlled_2006}):
\begin{equation}\label{controlled generator_pde}
\begin{aligned}
   \frac{\partial }{\partial
   t}\mathbb{E}^{x}(\phi(X_{t})) & =\mathcal{L}_{u}\mathbb{E}^{x}(\phi(X_{t})), &
   \mathbb{E}^{x}(\phi(X_{0})) &= \phi(x).
\end{aligned}
\end{equation}
Note that~\eqref{controlled generator_pde} is now a time-inhomogeneous equation.

%% file: Contents/2_Learning_methods.tex
To perform numerical calculations on the Koopman generator, one can compute a Galerkin projection onto a finite dimensional space $\mathcal{F}=\text{span}\{\psi_{i}\}_{i=1}^{N}$, spanned by a finite basis of observable functions. As shown in~[\cite{klus_data-driven_2020}], the matrix representation of this Galerkin projection is
$$
L=AC^{\dagger},
$$
with the stiffness and mass matrices 
$$A_{ij}=\langle\mathcal{L}\psi_{i},\psi_{j}\rangle_{\mu},\,\,\,\,C_{ij}=\langle\psi_{i},\psi_{j}\rangle_{\mu}.$$ Here, $\mu$ is a probability measure, which is typically chosen as an invariant measure in the time-homogeneous case. Moreover, if we have a data set $X = \{x_{l}\}_{l=1}^{m}$ of samples from $\mu$, we can compute empirical estimators for the stiffness and mass matrices. 

To this end, we introduce vector valued functions $\Psi\colon\mathbb{R}^{n}\to\mathbb{C}^{N}$ and $\mathcal{L}\Psi\colon\mathbb{R}^{n}\to\mathbb{C}^{N}$ as
$$
\begin{aligned}
        \Psi(\mathbf{x}) & = \begin{pmatrix}
            \psi_1(\mathbf{x}) \\  \vdots \\ \psi_N(\mathbf{x})
        \end{pmatrix}, &
        \cL\Psi(\mathbf{x}) &= \begin{pmatrix}
            \cL\psi_1(\mathbf{x}) \\  \vdots \\ \cL\psi_N(\mathbf{x})
        \end{pmatrix}.
\end{aligned}
$$
Likewise, we make the definitions
\begin{equation*}
    \begin{split}
            \Psi(X) &=\begin{bmatrix}
    \Psi(x_{1}) & \cdots & \Psi(x_{m})
\end{bmatrix} \in \mathbb{C}^{N \times m}, \\
        \mathcal{L}\Psi(X) &=\begin{bmatrix}
    \mathcal{L}\Psi(x_{1}) & \cdots & \mathcal{L}\Psi(x_{m})
\end{bmatrix} \in \mathbb{C}^{N \times m}.
    \end{split}
\end{equation*}

Then the empirical estimators for the mass and stiffness matrix are given by
$$
\hat{A}=\frac{1}{m} \cL\Psi(X)\Psi(X)^\htop,\,\,\,\,\hat{C}=\frac{1}{m}\Psi(X)\Psi(X)^\htop.
$$
The final empirical matrix approximation of Koopman generator is given by
$$
\hat{L}=\hat{A}\,\left(\hat{C} + \lambda\mathrm{Id}\right)^{\dagger},
$$
where $\lambda \geq 0$ is a regularization parameter. This estimation method is called generator extended dynamic mode decomposition (gEDMD)~[\cite{klus_data-driven_2020}].

\subsubsection{Prediction with gEDMD} 

Furthermore, if we assume that $\mathbb{E}^{x}(\phi(X_{t}))\in\mathcal{F}$ for $\phi \in \mathcal{F}$, i.e.
\begin{equation}\label{finite-dim exp}
\mathbb{E}^{x}(\phi(X_{t}))=\sum_{i=1}^{N}v_{i}(t)\psi_{i}(x)=V(t)^{\mathrm{H}}\Psi(x)
\end{equation}
with time-dependent vector $V(t)=(v_{1}(t),\cdots,v_{N}(t))^\top$, then (\ref{generator pde}) turns into
$$
\frac{\partial}{\partial t}V(t)^{\mathrm{H}}\Psi(x)=V(t)^{\mathrm{H}}\mathcal{L}\Psi(x),
$$
which implies that the time evolution of the expectation satisfies a data-driven linear ODE on the finite-dimensional space $\mathcal{F}$, i.e.
\begin{equation}\label{ode}
    \frac{d}{dt}V(t)^{\mathrm{H}}=V(t)^{\mathrm{H}}\hat{L}.
\end{equation}

\subsubsection{Bi-linear ODE for Control-affine Systems}
For a time-dependent system such as~(\ref{control stocha system}), constructing the empirical generator approximation becomes quite challenging because, according to~\eqref{controlled generator_pde}, $\mathcal{L}_{u}\psi_{k}(x)$ 
depends not only on the data but also on the time-dependent input signal $u(t)$. If we consider constant inputs 
$$
u\in\mathcal{U}=\{0,e_{1},\cdots,e_{p}\},
$$
where $e_{i}\in\mathbb{R}^{p}$ for $i=1,\cdots,p$ are unit vectors, then the generators corresponding to these constant inputs are given by
\begin{align*}
  \mathcal{L}_{0} & =b\nabla_{x}+\frac{1}{2}\Sigma\colon\nabla_{x}^{2}, &
  \mathcal{L}_{e_{i}} & =(b+G_{i})\nabla_{x}+\frac{1}{2}\Sigma\colon\nabla_{x}^{2}
\end{align*}
which, based on (\ref{controlled generator}), results in 
\begin{equation}\label{bilin generator}     \mathcal{L}_{u}= \mathcal{L}_{0}+\sum_{i=1}^{p}(\mathcal{L}_{e_{i}}-\mathcal{L}_{0})u_{i}(t).
\end{equation}
That is, the control-affine drift leads to a control-affine Koopman generator~[\cite{peitz_koopman_2019}]. Consequently, the evolution equation~(\ref{controlled generator_pde}) can be written as follows
\begin{equation}\label{generator_u_pde}
\left\{
\begin{aligned}
   &\frac{\partial }{\partial
   t}\mathbb{E}^{x}(\phi(X_{t}))=(\mathcal{L}_{0}+\sum_{i=1}^{p}(\mathcal{L}_{e_{i}}-\mathcal{L}_{0})u_{i}(t))\mathbb{E}^{x}(\phi(X_{t})),\\
   &\mathbb{E}^{x}(\phi(X_{0}))=\phi(x).
\end{aligned}\right.
\end{equation}
For given data set $X$, vector valued function $\Psi$ and input signal $u(t)$, the empirical matrix approximations of generators with constant inputs $L_{0}$ and $L_{e_{i}}$ for $i=1,\cdots,p$ can be obtained by applying gEDMD as outlined in~Alg.\ref{alg bilin generator}. Hence, by~(\ref{bilin generator}), the matrix approximation of the control-affine generator is likewise given by
$$
L_{u}=L_{0}+\sum_{i=1}^{p}(L_{e_{i}}-L_{0})u_i.
$$
Under the assumption (\ref{finite-dim exp}), and by a similar derivation as for (\ref{ode}), the PDE (\ref{generator_u_pde}) can be approximated by the following bi-linear ODE
\begin{equation}\label{bilin ode}
    \frac{d}{dt}V(t)^{\mathrm{H}}=V(t)^{\mathrm{H}}L_{u}.
\end{equation}

\begin{algorithm}[h]
\caption{gEDMD for Control-Affine Systems}\label{alg bilin generator}
\begin{algorithmic}[1]
\Require $\begin{aligned}
    X&=(x_1, x_2,\cdots,x_m)\in\mathbb{R}^{n\times m}\\
    \Psi&=(\psi_1, \cdots,\psi_N)^\top, \,
    u(\cdot)\in\mathbb{R}^{p},\, \lambda \geq 0
\end{aligned}$
\Ensure $L_{0},L_{e_{i}},L_{u}\in\mathbb{C}^{N\times N}$ for $i=1,\cdots,p$
\Statex
\State $C=\Psi(X)\Psi(X)^{\mathrm{H}}$
\State Set $G_0 = 0$ and $e_0 = 0$.
\For{$i=0,\cdots,p$}
\State $A_{e_{i}} = \cL_{e_i}\Psi(X)\Psi(X)^{\mathrm{H}}$
\State $L_{e_{i}}=A_{e_{i}} \left(C + \lambda \mathrm{Id}\right)^{\dagger}$
\EndFor
\State $
L_{u}=L_{0}+\sum_{i=1}^{p}(L_{e_{i}}-L_{0})u_i
$
\end{algorithmic}
\end{algorithm}

\subsubsection{Random Fourier Features} Numerous options exist for choice of basis functions $\{\psi_{i}\}_{i=1}^{N}$, e.g. Hermite polynomials and radial basis functions (RBFs) [\cite{williams_data-driven_2015}]. In order to generate expressive basis sets while keeping the number of parameters manageable, kernel-based methods relying on reproducing kernel Hilbert spaces (RKHS), see e.g. [Section 4 of~\cite{christmann_support_2008}], have been proposed as discussed in~[\cite{williams_kernel-based_2015}, \cite{klus_kernel-based_2018}]. The random Fourier features (RFFs) outlined in~[\cite{rahimi_random_2007}] can then be used to efficiently approximate these kernel functions.

We consider a translation invariant kernel $k$ on $\mathbb{X} = \mathbb{R}^n$, which satisfies $k(x,x) = 1$. Then, Bochner's theorem guarantees that the kernel is indeed the Fourier transformation of a finite positive Borel probability measure $\rho$ on $\mathbb{R}^{n}$, i.e.
\begin{equation*}
    \begin{split}
        k(x,y) &=\mathbb{E}_{\omega\sim\rho}(e^{-ix^\top\omega}\overline{e^{-iy^\top\omega}}).        
    \end{split}
\end{equation*}
To approximate the kernel, we define the random Fourier feature (RFF) mapping $\phi_{\text{RFF}}\colon\mathbb{R}^{n}\to\mathbb{C}^{N}$ by
\begin{equation}\label{rff}
\phi_{\text{RFF}}(x)=\begin{pmatrix}
    \phi_{1}(x)\\ \vdots\\ \phi_{N}(x)
\end{pmatrix}=\begin{pmatrix}
    e^{i x^\top\omega_{1}}\\ \vdots\\
    e^{i x^\top\omega_{N}}
\end{pmatrix},
\end{equation}
where $\{\omega_{i}\}_{i=1}^{N}$ are sampled from the spectral measure $\rho$. The point-wise approximation to the kernel $k$ is then given as $k(x, y) = \frac{1}{p}\phi_{\text{RFF}}(x)^\htop\phi_{\text{RFF}}(y)$. As explained in~[\cite{nuske_efficient_2023}], approximating kernelized gEDMD by random Fourier features is equivalent to just using standard gEDMD with the randomly generated basis set $\phi_{\text{RFF}}$. That is the approach we follow in this paper.

%% file: Contents/3_Solution_of_Control_Problem.tex
Our overall goal in this work is to find an optimal input $u^{\star}$ on the time horizon $[0,T]$ such that it solves a nonlinear optimal control problem (OCP) of the form:
\begin{equation}\label{cost_u}
\begin{aligned}
\min_{u:\,[0, T] \mapsto \mathbb{R}^{p}}  &\quad&  &J(x,u),\\
\text{s.t.} &\quad&  &\text{the dynamics }\eqref{control stocha system},\\
            &\quad& \text{with } &X_{0}=x.
\end{aligned}
\end{equation}
The cost functional we consider is typically of the form:
$$
    J(x,u):=\mathbb{E}^{x}\{\int_{0}^{T}l(X_{t},u(t))dt+\mathcal{T}(X_{T})\},
$$
where $l(\cdot,\cdot)$ and $\mathcal{T}(\cdot)$ are the running cost and the terminal cost, respectively. We only consider separable running costs of the form
$$
l(X_{t},u(t)) = l_{1}(X_{t}) + l_{2}(u(t)).
$$
and we approximate the time integral  using the piecewise trapezoidal rule.
The expectations can then be approximated using gEDMD, if we treat the expectations as the solutions to the PDE~(\ref{generator_u_pde}) with the observables $l_{1}(\cdot)$ and $\mathcal{T}(\cdot)$. To this end, the observables $l_1, \, \mathcal{T}$ must be contained in $\mathcal{F}$, i.e. there are coefficient vectors such that
\begin{align*}
    l_1(x) &= V_{l_1}^{\mathrm{H}} \Psi(x), & \mathcal{T}(x) &= V_{\mathcal{T}}^{\mathrm{H}} \Psi(x).
\end{align*}
Using Alg.~\ref{alg bilin generator} to provide a matrix approximation of the control-affine generator $L_{u}$, and based on the assumption (\ref{finite-dim exp}), the approximated expectations $\mathbb{E}^{x}(l_{1}(X_{t}))$ and $\mathbb{E}^{x}(\mathcal{T}(X_{t}))$ can be computed by solving the bi-linear ODE~\eqref{bilin ode} with the initial conditions $V_{l_{1}}(0) = V_{l_1}$ and $V_{\mathcal{T}}(0) = V_{\mathcal{T}}$, respectively.
We note that if we are only looking to solve~\eqref{cost_u} for a single initial position $x$, it is more efficient to solve a dual problem instead. Denoting the flow map of the ODE~\eqref{bilin ode} by $\Phi^t_u$, then for any any observable $\phi \in \mathcal{F}$, we have
\begin{equation*}
    \mathbb{E}^x[\phi(X_t)] = \left(\Phi^t_u (V_\phi(0)\right)^\htop \Psi(x).
\end{equation*}
For fixed $u$, the ODE~\eqref{bilin ode} is a linear time-dependent equation, and hence the flow map $\Phi^t_u$ is also linear. It follows that
\begin{align*}
    \mathbb{E}^x[\phi(X_t)] &= V_\phi(0)^\htop (\Phi^t_u)^\htop \Psi(x) = V_\phi(0)^\htop \gamma(t),
\end{align*}
where $\gamma(t) := (\Phi^t_u)^\htop \Psi(x)$ is the solution to the adjoint equation
$$
 \frac{d}{dt}\gamma(t)^\htop = \gamma(t)^\htop L_{u}^{\htop},
$$
with the initial condition $\gamma(0)^\htop = \Psi(x)^\htop$. Putting these pieces together, we can evaluate the cost functional $J(x, u)$ for any input signal $u$. The optimization is then carried out using a black-box solver from the \emph{scipy} library.

%% file: Contents/4_Results.tex
In this section, we perform numerical simulations using the following one-dimensional SDE with control-affine drift, given by
\begin{equation}\label{dw_sde}
    dX_{t}=-(\nabla V(X_{t})+K_{\text{bias}}(X_{t}-u(t)))dt+\sqrt{2\beta^{-1}}
    dW_{t}
\end{equation}
with the double-well potential $V\colon\mathbb{R}\to\mathbb{R}$ defined as
\begin{equation}\label{dw_potential}
V\colon x\mapsto K_{\text{dw}}(x^{2}-1)^{2}.
\end{equation}
Here, $K_{\text{dw}},\, K_{\text{bias}},\,\beta \in\mathbb{R}_{>0}$, are constants, and $u\colon\,[0, T]\to\mathbb{R}$ is the input signal. The diffusion constant $\beta$ is called inverse temperature in statistical physics, we set $\beta = 1$ in all numerical experiments. For training, we use the constant inputs $u \in \{-1,1\}$. Note that the corresponding drift fields are the derivatives of quadratic bias energies
$$
B_{1}:=\frac{K_{\text{bias}}}{2}(x+1)^{2},\,\,\,\,B_{2}:=\frac{K_{\text{bias}}}{2}(x-1)^{2},
$$
which are used in many adaptive sampling algorithms. Also, note that we do not learn a model for the generator $\cL_0$ of the metastable physical system, but use two dynamics with localized energy functions instead.
\begin{figure}[t]
    \centering
    \subfigure[]{\includegraphics[width=0.48\linewidth]{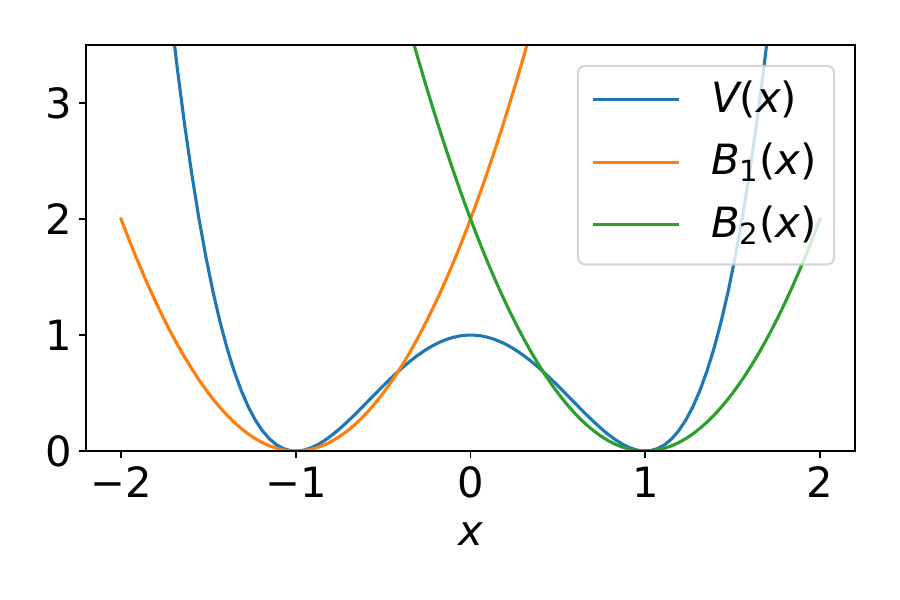}}
    \subfigure[]{\includegraphics[width=0.48\linewidth]{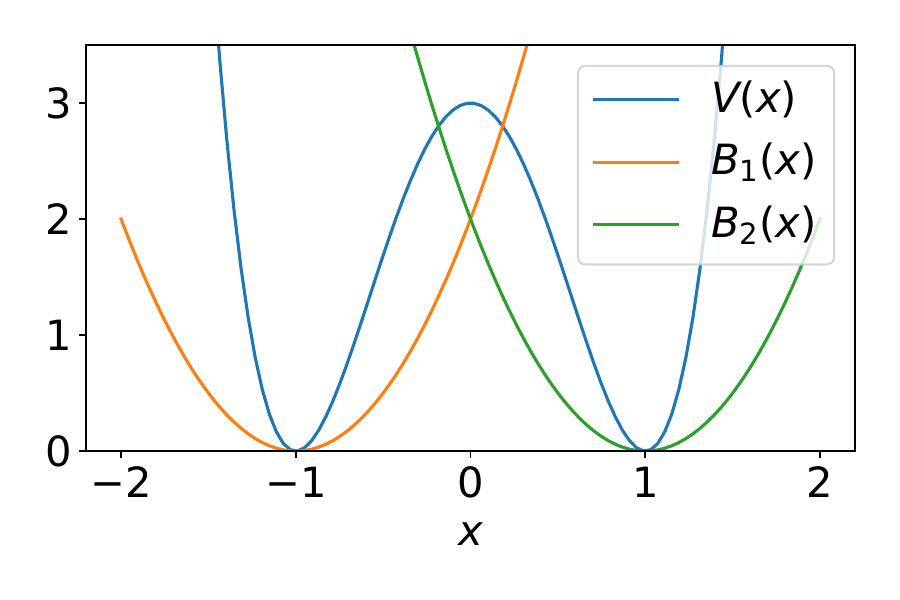}}
    \caption{Double-well potential $V$ and bias energies $B_i$ with $K_{\text{bias}}=4$ for  $K_{\text{dw}}=1$ (a) and $K_{\text{dw}}=3$ (b)}
    \label{fig:potentials}
\end{figure}

%% file: Contents/5_Prediction.tex
As a first application of the bi-linear gEDMD framework, we consider a fixed input $u(\cdot)$, and use gEDMD to predict the expectation of observable functions along the trajectory of the dynamics. For the double-well example, we choose the input as a periodic signal $u(t):=\cos(2t)$. We set
the initial condition $x=0.5$, the step size $\triangle t=10^{-3}$ and the
number of steps $N_{s}=5\times 10^{3}$. Additionally, we use random Fourier
features $\phi_{\text{RFF}}$ (\ref{rff}) to span the finite-dimensional
subspace, where $\{\omega_{i}\}_{i=1}^{N}$ are sampled from the spectral density
of a Gaussian kernel with the bandwidth $0.5$, and the number of features is set
to $N=50$. Furthermore, the training
data $\{x_{l}\}_{l=1}^{m}$ are $m$ points uniformly distributed on the interval
$[-2,2]$. On these data, we learn the matrix models for the generators
$\mathcal{L}_{-1}$ and $\mathcal{L}_{1}$ for the constant inputs $u=-1$ and $u=1$, respectively.

Using the bi-linear dynamics~\eqref{bilin ode}, we compute an approximation $\hat{\mathbb{E}}^{x}(X_{t})$ to the true expectation $\mathbb{E}^{x}(X_{t})$ of the dynamics~(\ref{dw_sde}). As a ground truth for $\mathbb{E}^{x}(X_{t})$, we compute empirical averages over 100 trajectories, generated by the Euler-Maruyama scheme. We denote the absolute value of the prediction error by
$$
|e(t)|:=|\mathbb{E}^{x}(X_{t})-\hat{\mathbb{E}}^{x}(X_{t})|.
$$
We investigate how the constants $K_{\text{dw}}$, $K_{\text{bias}}$, the regularization parameter $\lambda$ and the size of the training data $m$ affect the absolute error $|e(t)|$. In addition, we consider the prediction as a failure if $|e(t)| \ge 1$ for any $t$. We omit failed predictions when computing mean errors for $e(t)$, but we also report the success rate $\delta$.
\begin{figure}[t]
    \centering
        \subfigure[]{\includegraphics[width=0.48\linewidth]{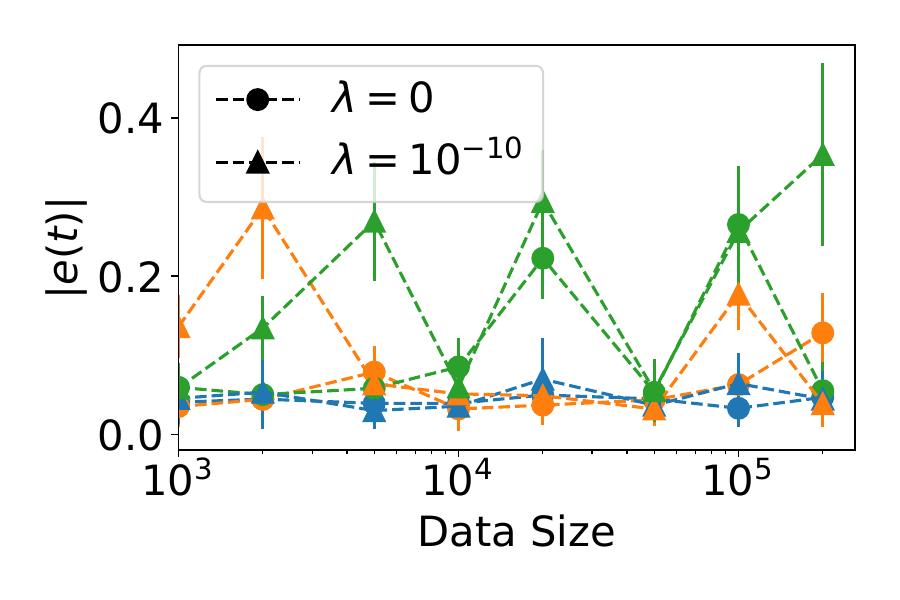}}
        \subfigure[]
        {\includegraphics[width=0.48\linewidth]{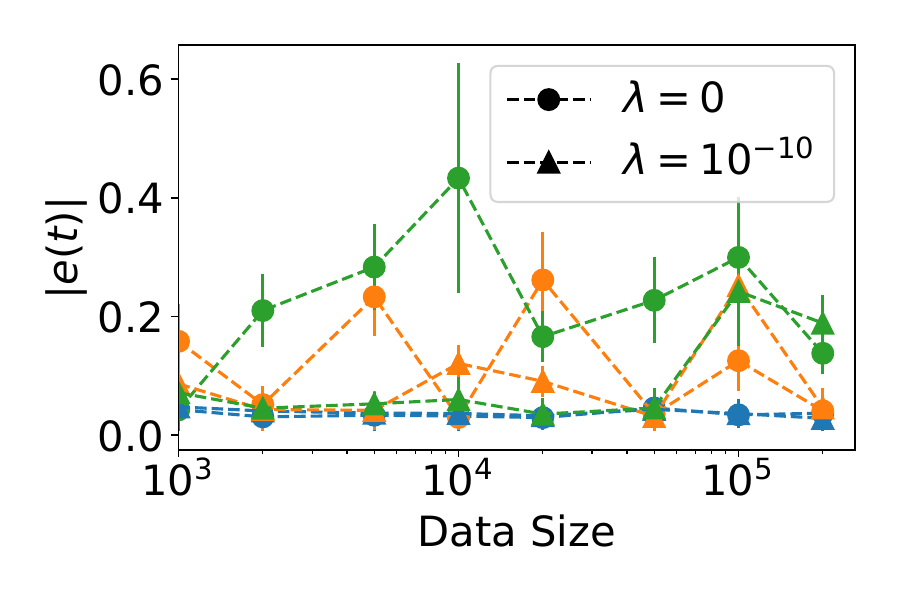}}
        \subfigure[]
        {\includegraphics[width=0.48\linewidth]{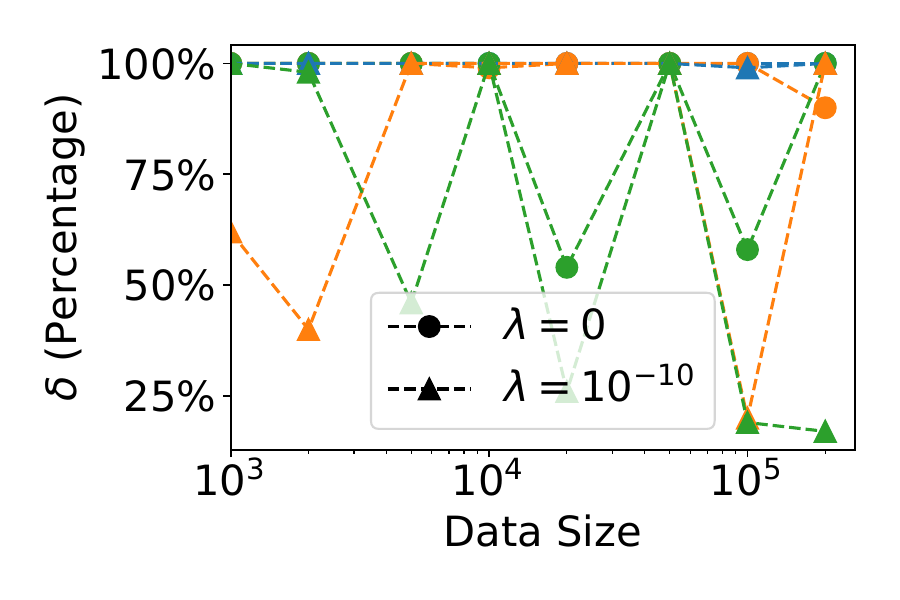}}
        \subfigure[]
        {\includegraphics[width=0.48\linewidth]{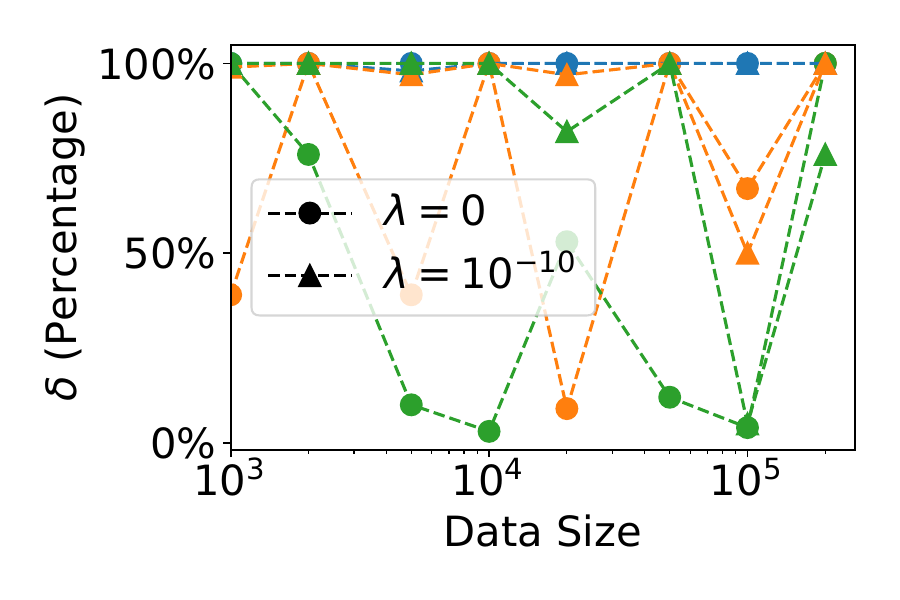}}
    \caption{Absolute prediction errors (top row) and success rates $\delta$ (bottom row) with $K_{\text{bias}}=3$~(left) and $K_{\text{bias}}=4$~(right). Blue, orange and green lines stand for the cases when $K_{dw}=\{1,2,3\}$ respectively.}
    \label{fig:prediction_error}
\end{figure}


Fig.~\ref{fig:prediction_error} shows how the error $|e(t)|$ (panels a\&b) and the success rate (panels c\&d) vary with the data size for each value of $K_{\text{dw}}=\{1,2,3\}$, regularization parameters $\lambda \in \{0,10^{-10}\}$, while keeping either the biasing constant at $K_{\text{bias}} = 3$ or at $K_{\text{bias}} = 4$. We observe that the predictions are very reliable for $K_{\text{dw}} \in \{1, 2\}$, and more susceptible to instabilities for the more metastable system $K_{\text{dw}} = 3$. We notice that for this latter case, using no regularization $\lambda = 0$ works better for $K_{\text{bias}} = 3$, while a small amount of regularization $\lambda = 10^{-10}$ seems appropriate for $K_{\text{bias}} = 4$.


 As a consequence, we will use the settings $(K_{\text{bias}},\lambda) = (3, 0)$ and $(K_{\text{bias}},\lambda) = (4, 10^{-10})$, and only moderate training data sizes around $m = 10^3$, in all of the following examples. A more detailed study of the interplay between basis set size $N$, data size $m$, and regularization parameter $\lambda$ is planned for the future.

%% file: Contents/6_Tracking.tex
We now consider a more complex problem: steering the expectation to a predefined value. This can be formulated as an optimization problem. Given a reference signal $x_{\text{ref}}\colon\,[0, T]\to\mathbb{R}$, the goal is to find an optimal input control $u^{\star}$ that solves the following finite horizon optimization problem:
\begin{equation}\label{track opt}
\begin{aligned}
\min_{u:\,[0, T] \mapsto \mathbb{R}}  &\quad&  &\int_{0}^{T}||\mathbb{E}^{x}(X_{t})-x_{\text{ref}}(t)||^{2}dt\\
\text{s.t.} &\quad&  &\text{the dynamics }\eqref{dw_sde},\\
            &\quad& \text{with } &X_{0}=x.
\end{aligned}
\end{equation}
We define the absolute tracking error as
$$
|e_{t}(t)|:=|\mathbb{E}^{x}(X_{t})-x_{\text{ref}}(t)|,
$$
while the reference signal is $x_{\text{ref}}(t)=\cos(2t)$, i.e. the same
signal that was used as fixed input before. We consider four different settings
for the system parameter $K_{\text{dw}}$, the biasing strength
$K_{\text{bias}}$, and the regularization parameter $\lambda$, namely
$(K_{\text{dw}},K_{\text{bias}},\lambda) \in \{(1,3,0),(3,3,0),(1,4,10^{-10}),(3,4,10^{-10})\}$. We set the data size to $m=10^{3}$ and the final time to $T = 2$, thus tracking the reference signal for one complete transition from the left to the right minimum. The remaining parameter values are the same as those specified in the previous section.

\begin{figure}[t]
    \centering
    \begin{minipage}[t]{0.48\linewidth}
        \centering
    \subfigure[]{\includegraphics[width=\linewidth]{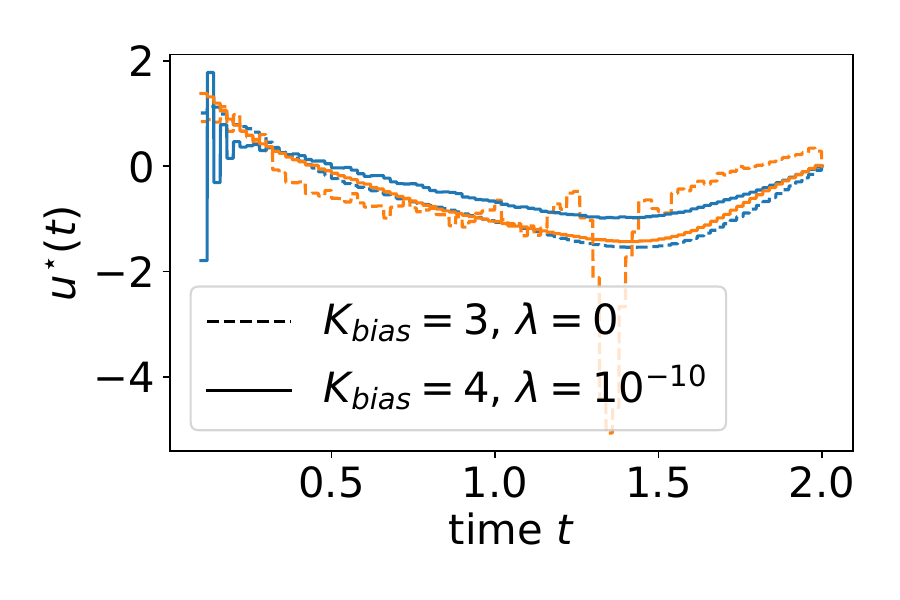}}
    \subfigure[]{\includegraphics[width=\linewidth]{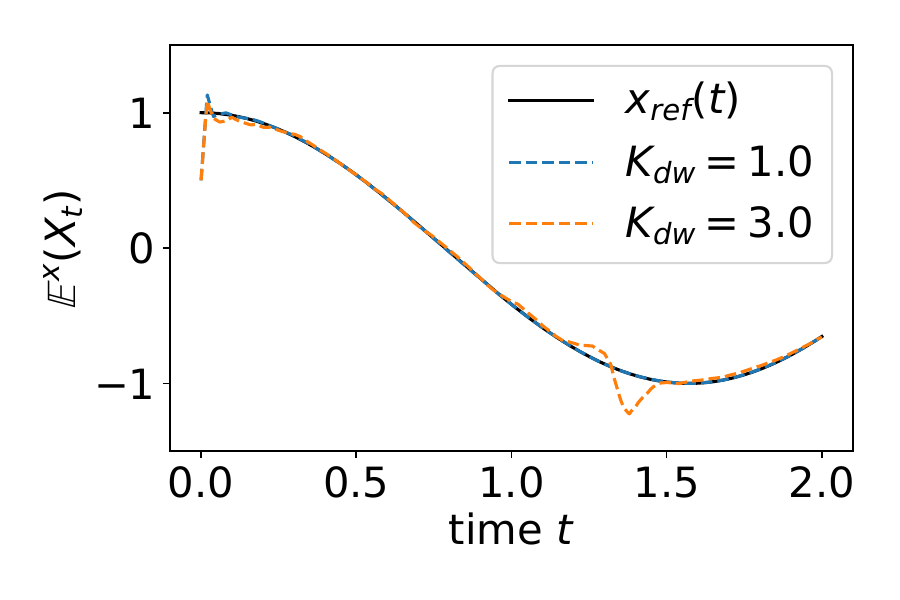}}
    
        \end{minipage}
    \hfill
    \begin{minipage}[t]{0.48\linewidth}
        \centering
    \subfigure[]{\includegraphics[width=\linewidth]{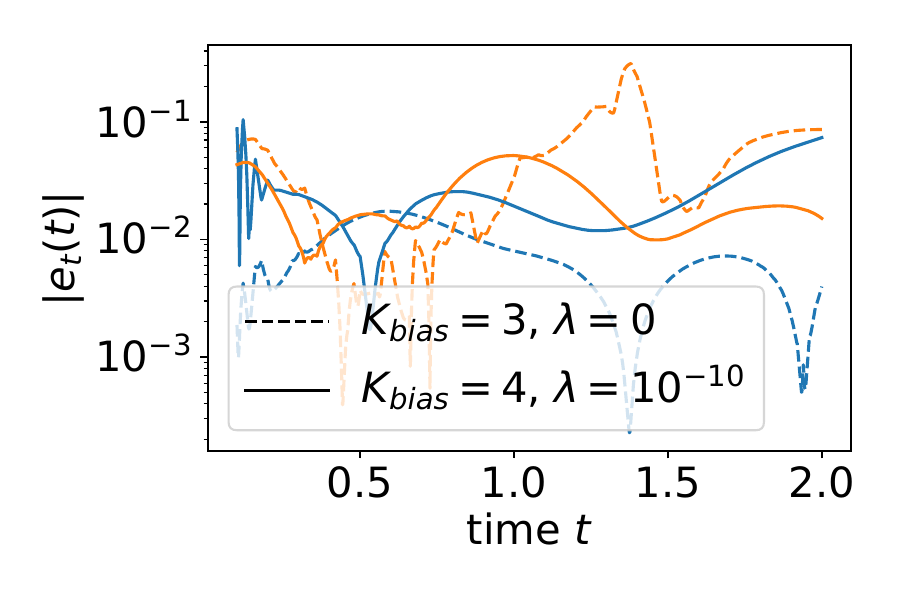}}
    \subfigure[]{\includegraphics[width=\linewidth]{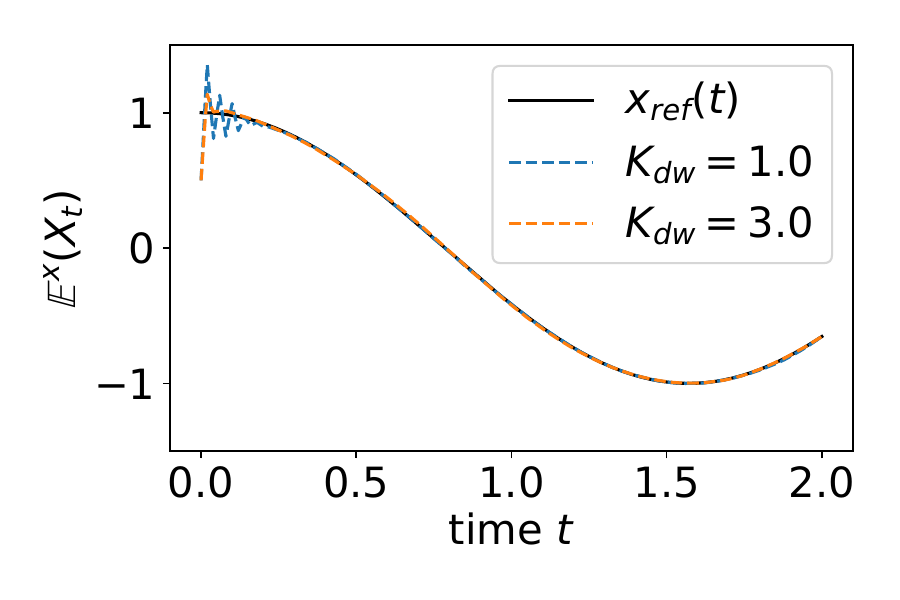}}

    \end{minipage}
    \caption{Numerical results for the tracking optimal control problem~(\ref{track opt}). (a)\&(c) are the optimal signals and tracking errors for four parameter settings.  Blue and orange lines stand for the systems when $K_{dw} \in \{1,3\}$, respectively. Panels (b)\&(d) are the tracking performances for fixed $K_{\text{bias}}=3$ and $K_{\text{bias}}=4$, respectively.}
    \label{fig:tracking}
\end{figure}

For four settings, we show in the first row of Fig.~\ref{fig:tracking}, the piecewise continuous solution $u^{\star}$ and the absolute error $|e_{t}|$. The corresponding tracking performances are shown in the second row of Fig.~\ref{fig:tracking}. As the initial condition $x$ and $x_{\text{ref}}(0)$ are not identical, we omit the time interval $[0, 0.1]$, during which the controlled system catches up with the reference, in the first row.


The results demonstrate that $\mathbb{E}^{x}(X_{t})$ effectively tracks the reference $x_{\text{ref}}(t)$ over the time interval $[0, 2]$ in all cases. For the systems with $K_{\text{dw}} = 1$, accurate tracking is easily achieved. We only observe that for $K_{\text{bias}}=4,\lambda=10^{-10}$, the optimal solution displays some oscillations around the reference initially, while for $K_{\text{bias}}=3, \lambda=0$, the tracking error is limited to less than one per cent for most of the control horizon. For the second, more strongly metastable system with $K_{\text{dw}} = 3$, the solution for $K_{\text{bias}}=3,\lambda=0$ allows for some stronger oscillation, while the tracking error of the system with $K_{\text{bias}}=4,\lambda=10^{-10}$ is limited to less than 10 per cent throughout.

%% file: Contents/7_Loss_Functions.tex
In this section, we explore the use of optimal control to accelerate transitions between metastable states. For the double-well potential~(\ref{dw_potential}), which has two metastable minima at $x = -1$ and $x = 1$, we consider the problem of enforcing a transition from the left minimum to the right basin. That is, with the initial condition $x = -1$, we aim to solve the optimal control problem~\eqref{cost_u} with the cost defined as
$$
J(x,u):=\mathbb{E}^{x}\{\int_{0}^{T}l(X_{t},u(t))dt\}+(1-\mathbb{E}^{x}\{X_{T}\})^{2}.
$$
The terminal cost is chosen to ensure that the dynamics will approximately reach the right minimum at the final time $T$. We explore, for $c\in\mathbb{R}_{>0}$, two choices of the running cost: first,
\begin{equation}\label{running dw}
    l_{\text{dw}} := V(X_{t}) + c||u(t)||^{2},
\end{equation}
which penalizes the potential energy and the magnitude of the input $u$. Another meaningful choice of running cost is
\begin{equation}\label{running bias}
    l_{\text{bias}}:=V(X_{t}) + c||X_{t}-u(t)||^{2},
\end{equation}
which exerts a penalty on both the physical energy and the bias energy spent during the simulation. We study four different system settings, namely $(K_{\text{dw}},K_{\text{bias}}) \in \{(1,3),(1,4),(3,3),(3,4)\}$. Furthermore, we fix the size of the training data as $m=10^{3}$ and the final time as $T = 1$. All remaining settings are the same as in section~\ref{sec prediction}.

%% file: Contents/8_Comparison.tex
\begin{figure}[t]
    \centering
    \begin{minipage}[t]{0.48\linewidth}
        \centering
        \text{Optimal Signals}
    \subfigure[$(1,3)$]{\includegraphics[width=\linewidth]{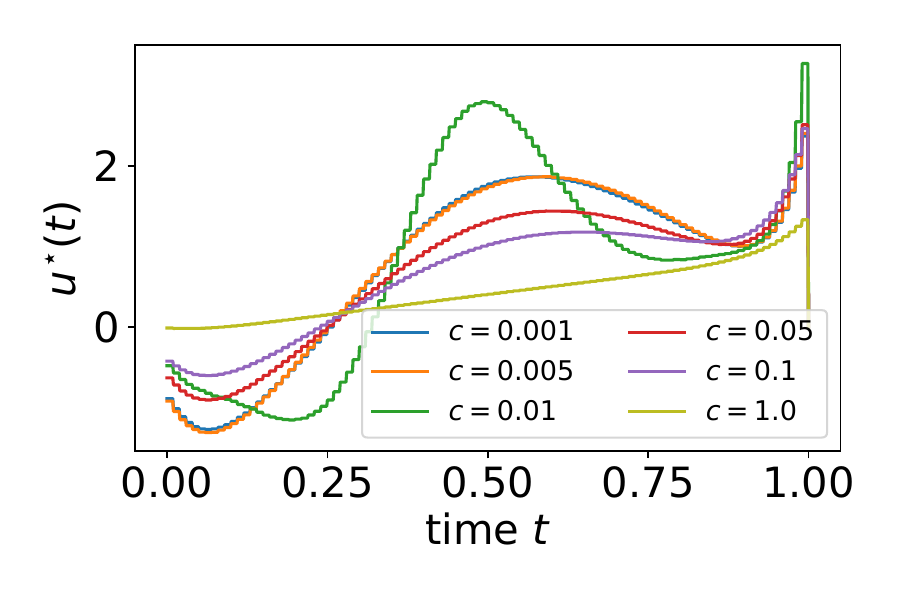}}
    \subfigure[$(3,3)$]{\includegraphics[width=\linewidth]{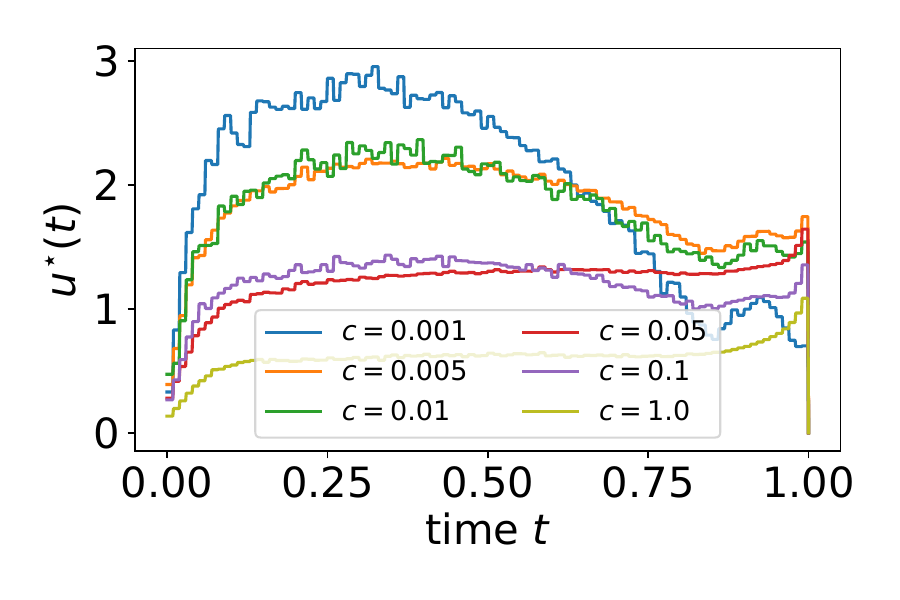}}
    \subfigure[$(3,4)$]{\includegraphics[width=\linewidth]{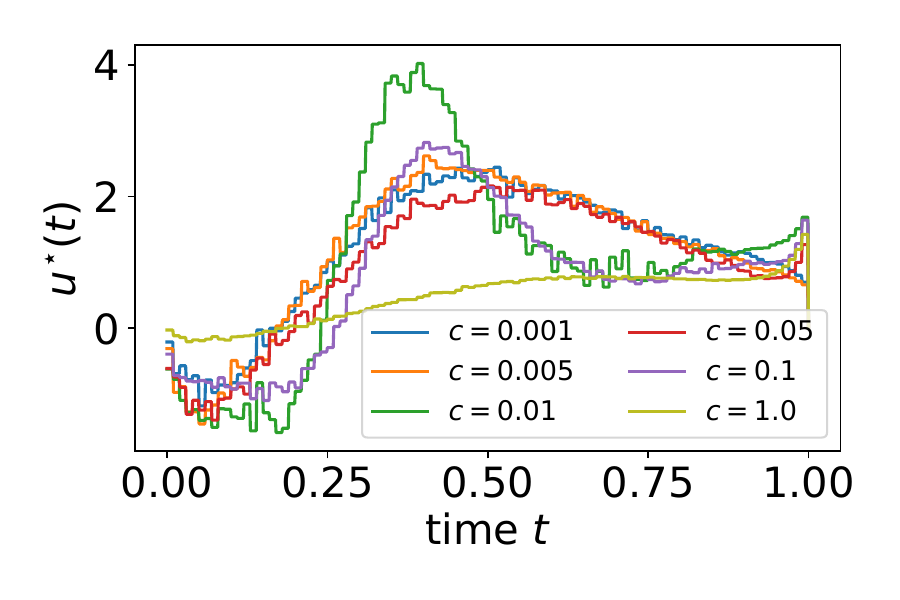}}
        \end{minipage}
    \hfill
    \begin{minipage}[t]{0.48\linewidth}
        \centering
        \text{Expectations}
    \subfigure[$(1,3)$]{\includegraphics[width=\linewidth]{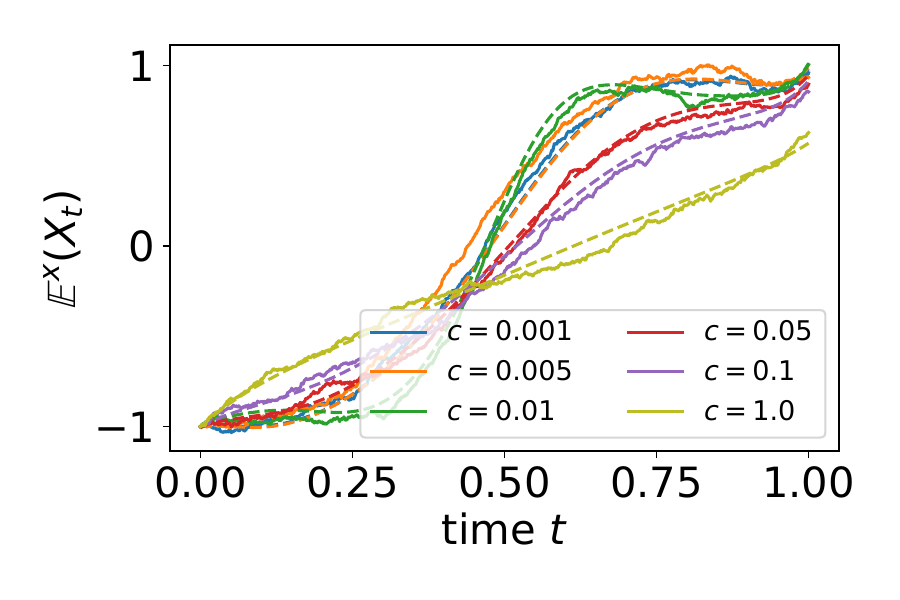}}
    \subfigure[$(3,3)$]{\includegraphics[width=\linewidth]{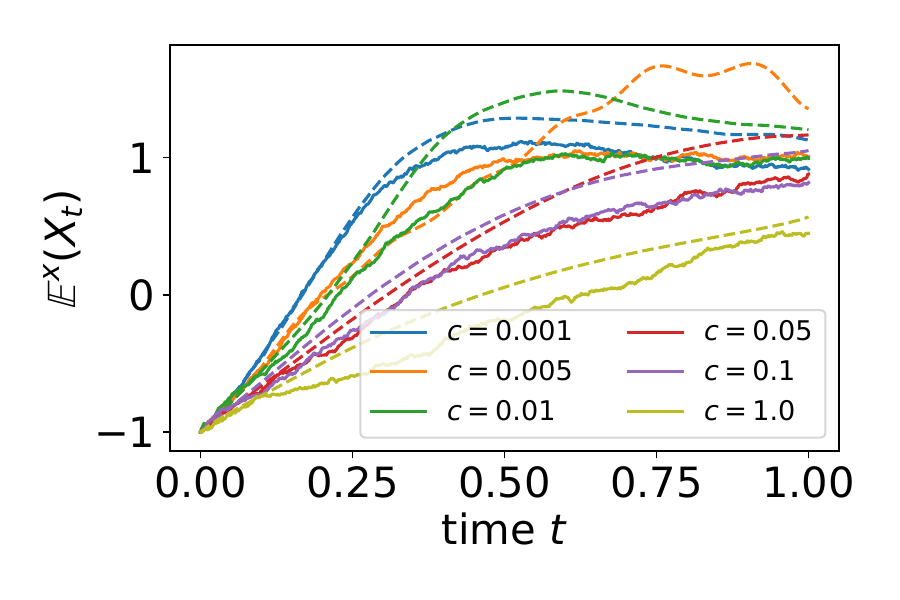}}
    \subfigure[$(3,4)$]{\includegraphics[width=\linewidth]{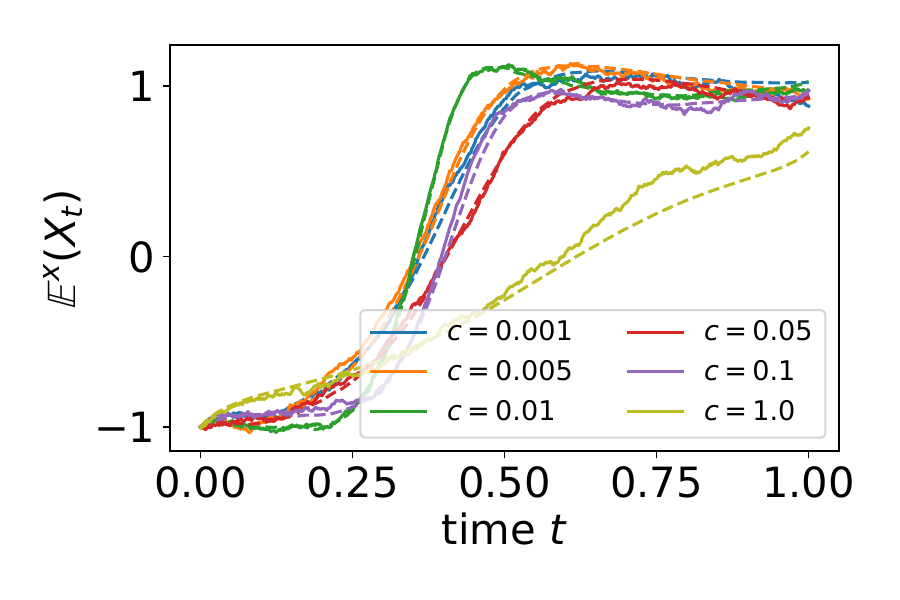}}
    \end{minipage}

    \caption{Numerical results of the problem (\ref{cost_u}) for the running cost (\ref{running dw}) with a range of parameters $c$. Left column are optimal signals and right column are expectations of the state $\mathbb{E}^{x}(X_t)$. From top to bottom are the results for the settings $(K_{\text{dw}},K_{\text{bias}}) \in \{(1,3),(3,3),(3,4)\}$, respectively, we omit the combination $(1, 4)$ as it is very similar to the first row. Solid and dashed lines stand for the reference and approximated expectations, respectively.}
    \label{fig:opt_dw}
\end{figure}
\begin{figure}[t]
    \centering
    \begin{minipage}[t]{0.48\linewidth}
        \centering
        \text{Optimal Signals}
    \subfigure[$(1,3)$]{\includegraphics[width=\linewidth]{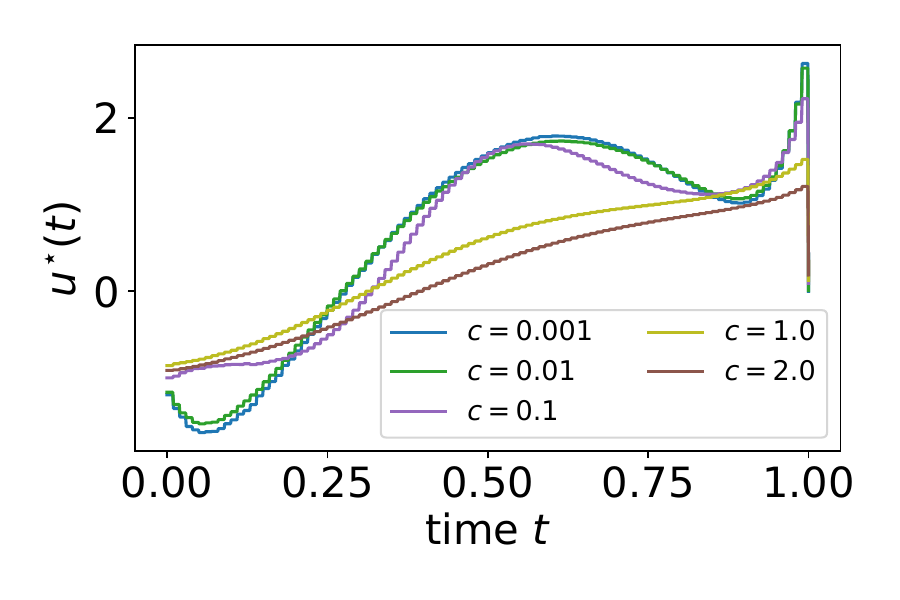}}
    \subfigure[$(3,3)$]{\includegraphics[width=\linewidth]{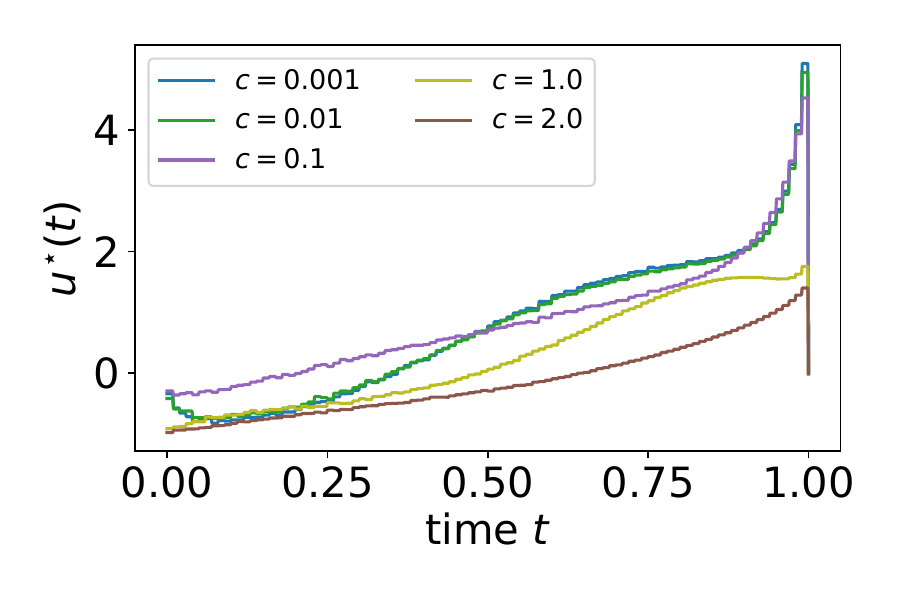}}
    \subfigure[$(3,4)$]{\includegraphics[width=\linewidth]{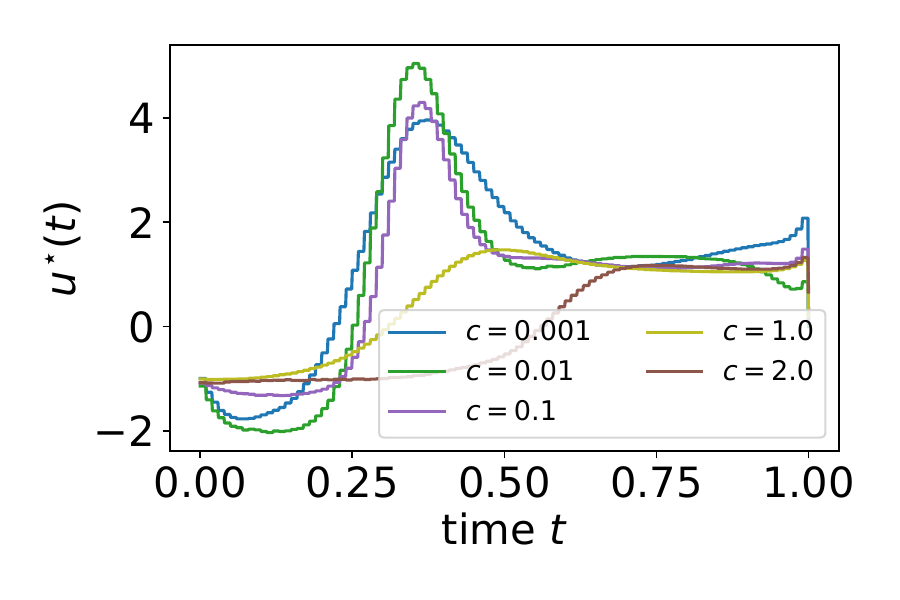}}
        \end{minipage}
    \hfill
    \begin{minipage}[t]{0.48\linewidth}
        \centering
        \text{Expectations}
    \subfigure[$(1,3)$]{\includegraphics[width=\linewidth]{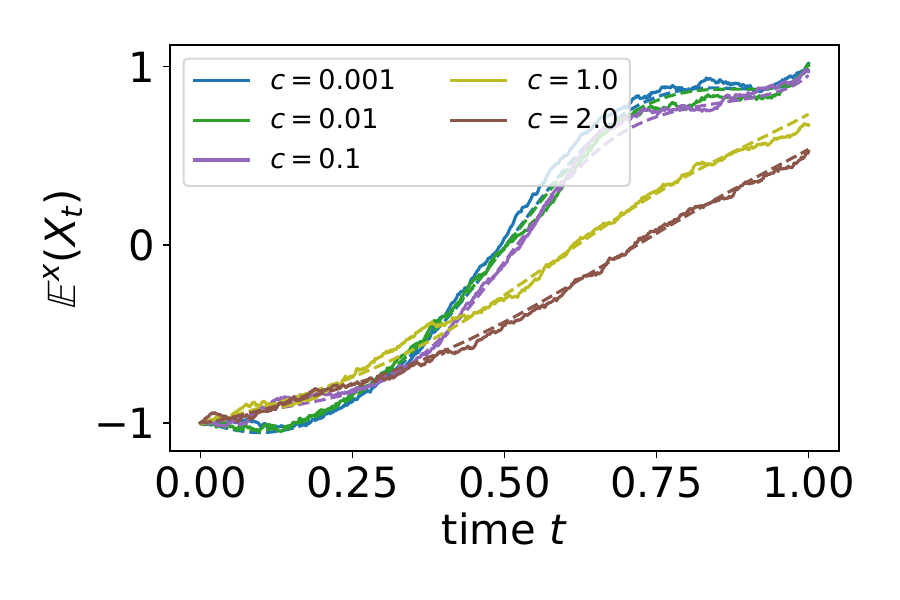}}
    \subfigure[$(3,3)$]{\includegraphics[width=\linewidth]{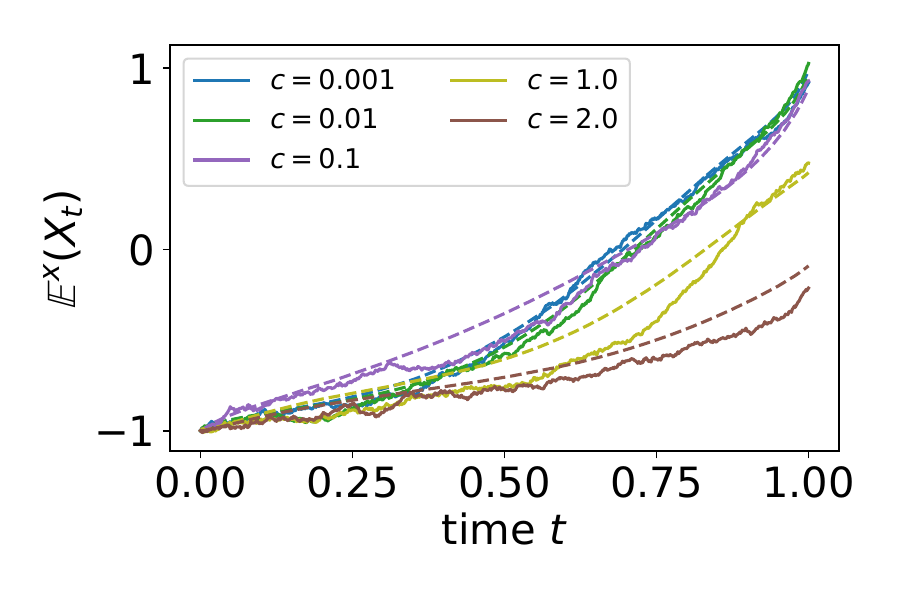}}
    \subfigure[$(3,4)$]{\includegraphics[width=\linewidth]{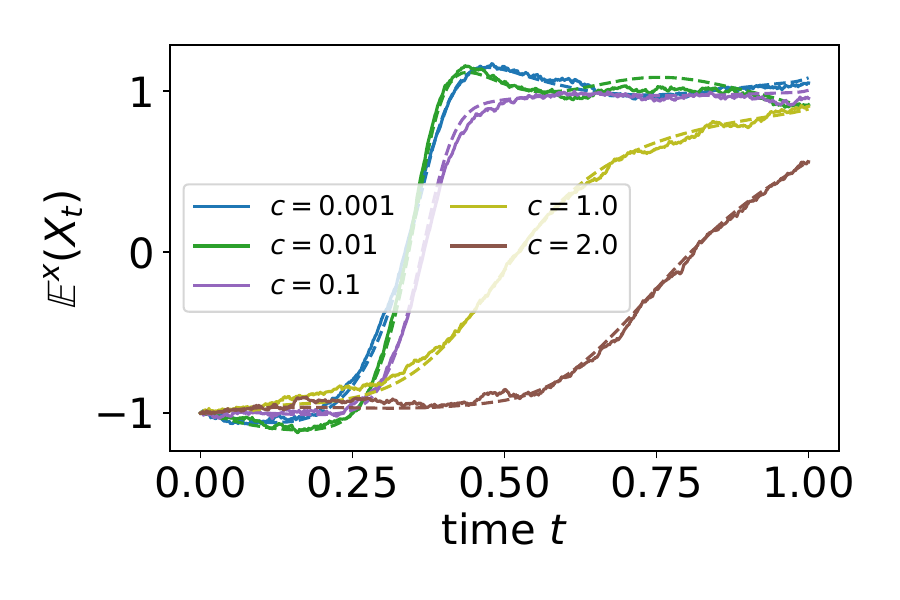}}
    \end{minipage}

    \caption{Same as in Figure~\ref{fig:opt_dw}, but for the running cost (\ref{running bias}).}
    \label{fig:opt_bias}
\end{figure}

For the four different system settings, we solve the optimization problems for a range of values of the parameter $c$ between $c = 10^{-3}$ and $c = 2.0$. The results are shown in Figs.~\ref{fig:opt_dw} \&~\ref{fig:opt_bias}. Generally speaking, the results for both loss functions and systems are fairly
similar. We find that there is a parameter range with $c$ between $10^{-3}$ and
$10^{-1}$, where the control objective of steering the system across the barrier is achieved. In most cases, the optimal signal displays a peak of the input signal, resulting in a rapid transition of the system state. The transition occurs a bit earlier for $(K_{\text{dw}},K_{\text{bias}}) = (3, 4)$ as opposed to $(K_{\text{dw}},K_{\text{bias}}) = (1, 3)$. For the parameters $(K_{\text{dw}},K_{\text{bias}}) = (3, 3)$, however, a more gradual transition is achieved by applying a rather flat input signal.
We also notice once again that in almost all cases, the expected values predicted by the bi-linear gEDMD model perfectly match the empirical mean values obtained from simulating the controlled dynamics~\eqref{control stocha system} with the optimal inputs signals $u^\star$.

%% file: Contents/9_Conclusion.tex
We have demonstrated the capabilities of the gEDMD algorithm for prediction and optimal control of control-affine stochastic systems. In particular, we have shown that optimal control policies for rare event sampling can be determined with limited amounts of data. Future research will explore the application of the method to more complex systems in higher dimensions, with a focus on the selection of training systems and generation of training data. We will also continue to investigate the efficient tuning of hyper-parameters, as well as the theoretical foundations of the method. We hope that our findings in this paper will inspire further research and practical applications in the field.